\documentclass[10pt]{article}
\usepackage[latin1]{inputenc}
\usepackage{amsmath, amsfonts, amssymb,amsthm}
\usepackage{todo}				
\usepackage{mathtools}
\usepackage{dsfont}

\usepackage{color}

\newtheorem{lemma}{Lemma}
\newtheorem{theorem}{Theorem}

\newtheorem{corollary}{Corollary}
\newtheorem{proposition}{Proposition}

\newcommand{\fp}{\mathbb{F}_p}
\renewcommand{\d}{\mathrm d}
\newcommand\Y{\mathbb {Y}}
\newcommand\ZZ{\mathbb{Z}}
\newcommand{\E}{\mathbb E}
\renewcommand\Pr{\mathbb{P}}
\title{Bounds for generalized Sidon sets}

\author{Xing Peng \thanks{Department of Mathematics, University of California, San Diego,  La Jolla, CA 92093, USA.
{\tt (x2peng@ucsd.edu)}. Research is supported in part by ONR MURI N000140810747, and AFSOR AF/SUB 552082. }
\and
Rafael Tesoro \thanks{Departamento de Matem\'{a}ticas, Universidad
Aut\'onoma de Madrid. 28049 Madrid, Spain. {\tt (rafael.tesoro@estudiante.uam.es)}}
\and
Craig Timmons
\thanks{Department of Mathematics, University of California, San Diego,  La Jolla, CA 92093, USA.
({\tt ctimmons@ucsd.edu}). Research is partially supported by NSF Grant DMS-1101489 through Jacques Verstra\"{e}te.
}
}
\date{}
\begin{document}
\maketitle

\begin{abstract}
Let $\Gamma$ be an abelian group and $g \geq h \geq 2$ be integers.  A set $A \subset \Gamma$ is a
$C_h[g]$-set if given any set $X \subset \Gamma$ with $|X| = k$, and any set $\{ k_1 , \dots , k_g \} \subset \Gamma$, at least one of the translates $X+ k_i$ is not contained in $A$.  For any $g \geq h \geq 2$, we prove that if $A \subset \{1,2, \dots ,n \}$ is a $C_h[g]$-set in $\mathbb{Z}$, then $|A| \leq (g-1)^{1/h} n^{1 - 1/h} + O(n^{1/2 - 1/2h})$.  We show that for any integer $n \geq 1$, there is a $C_3 [3]$-set $A \subset \{1,2, \dots , n \}$ with $|A| \geq (4^{-2/3} + o(1)) n^{2/3}$.
We also show that for any odd prime $p$, there is a $C_3[3]$-set $A \subset \mathbb{F}_p^3$ with
$|A| \geq p^2  - p$, which is asymptotically best possible.
Using the projective norm graphs from extremal graph theory, we show that for each integer $h \geq 3$, there is a $C_h[h! +1]$-set $A \subset \{1,2, \dots , n \}$
with $|A| \geq ( c_h +o(1))n^{1-1/h}$.
 A set $A$ is a \emph{weak $C_h[g]$-set} if we add the condition that the translates $X +k_1, \dots , X + k_g$ are all pairwise disjoint.
We use the probabilistic method to construct weak $C_h[g]$-sets in $\{1,2, \dots , n \}$ for any $g \geq h  \geq 2$. Lastly we obtain  upper bounds on infinite $C_h[g]$-sequences.  We prove that for any infinite $C_h[g$]-sequence $A \subset \mathbb{N}$, we have $A(n) = O ( n^{1 - 1/h} ( \log n )^{ - 1/h} )$ for infinitely many $n$, where $A(n) = | A \cap \{1,2, \dots , n \}|$.

\end{abstract}


\section{Introduction}
Given an integer $n \geq 1$, write $[n]$ for $\{1,2, \dots , n \}$.  Let $\Gamma$ be an abelian group and $g \geq h \geq 2$ be integers.  A set $A \subset \Gamma$ is a \emph{$C_h[g]$-set} if given any set $X \subset \Gamma$ with $|X| = k$, and any set $\{ k_1 , \dots , k_g \} \subset \Gamma$, at least one of the translates
\[
X+ k_i:=\{ x  + k_i : x \in X \}
\]
is not contained in $A$.
These sets were introduced by Erd\H os and Harzheim in \cite{ErdHar-1986}, and they are a natural generalization of the well-studied Sidon sets.  A Sidon set is the same as a $C_2[2]$-set.  We will always assume that $g \geq h \geq 2$.  The reason for this is that
if $X = \{ x_1 , \dots , x_k \}$ and $K = \{ k_1 , \dots , k_g \}$, then $A$ contains each of the translates
$X + k_1 , \dots , X + k_g$ if and only if $A$ contains each of the translates $K + x_1 , \dots , K + x_k$.

Our starting point is a connection between $C_h[g]$-sets and the famous Zarankiewicz problem from extremal combinatorics.  Given integers $m,n,s,t$ with $m \geq s \geq 1$ and $n \geq t \geq 1$, let $z(m,n,s,t)$ be the largest integer $N$ such that there is an $m \times n$ 0-1 matrix $M$, that contains $N$ 1's, and does not contain an $s \times t$ submatrix of all 1's.  Determining $z(m,n,s,t)$ is known as the problem of Zarankiewicz.

\begin{proposition} \label{prop1}
Let $\Gamma$ be a finite abelian group of order $n$.  Let $A \subset \Gamma$ and let $g \geq h \geq 2$ be integers.  If
$A$ is a $C_h[g]$-set in $\Gamma$, then
\begin{equation}\label{z eq}
n |A| \leq z(n,n,g,h).
\end{equation}
\end{proposition}
To see this, let $A \subset \Gamma$ be a $C_h[g]$-set where $\Gamma = \{ b_1 , \dots , b_n \}$ is a finite abelian group of order $n$.  Define an $n \times n$ 0-1 matrix $M$ by putting a 1 in the $(i,j)$-entry if
$b_i + b_j \in A$, and 0 otherwise.  A $g \times h$ submatrix of all 1's consists of a set
$X = \{x_1 , \dots , x_h \}$ of $h$ distinct elements of $\Gamma$, and a sequence $k_1 , \dots , k_g$ of $g$ distinct elements of $\Gamma$, such that $x_i + k_j \in A$ for all $1 \leq i \leq h$, $1 \leq j \leq g$.  There is no such submatrix since $A$ is a $C_h[g]$-set.  Furthermore, each row of $M$ contains $|A|$ 1's so that $n|A| \leq z(n,n,g,h)$.

F\"{u}redi \cite{fur1} proved that
\begin{equation}\label{fur ineq}
z(m,n,s,t) \leq ( s - t)^{1/2} n m^{1  - 1/t} + tm^{2 - 2/t} + tn
\end{equation}
for any integers $m \geq s \geq t \geq 1$ and $n \geq t$.  Therefore, if $A \subset \Gamma$ is a $C_h[g]$-set and $\Gamma$ is a finite abelian group of order $n$, then
\begin{equation}\label{group bound}
|A| \leq (g - h + 1)^{1/h} n^{1-1/h} + hn^{1 - 2/h} + h.
\end{equation}
If $A \subset [n]$ is a $C_h[g]$-set, then it is not difficult to show that $A$ is a $C_h[g]$-set in $\mathbb{Z}_{2n}$, thus by
(\ref{group bound}),
\[
|A| \leq (g - h + 1)^{1/h} 2^{1 - 1/h}n^{1 - 1/h} + h (2n)^{1 - 2/h} +h.
\]
Our first result improves this upper bound.

\begin{theorem}\label{thm1}
If $A \subset [n]$ is a $C_h[g]$-set with $g \ge h \ge 2$,  then
\begin{equation} \label{eq:BetterUpperBound-for-F_h(g,n)}
|A| \le (g-1)^{1/h} n^{1- 1/h} + O\left(n^{ 1/2 - 1/2h} \right).
\end{equation}
\end{theorem}

This theorem is a refinement of the estimate $|A|=O(n^{1-1/h})$ proved by Erd\H{o}s and Harzheim \cite{ErdHar-1986}.
Recall that $C_2[2]$-sets are Sidon sets.  Theorem \ref{thm1} recovers the well-known upper bound
for the size of Sidon sets in $[n]$ obtained by Erd\H os and Tur\'an \cite{ErdTur-1941}.
In general, $C_2[g]$-sets are those sets $A$ such that each nonzero difference $a-a'$ with $a ,a' \in A$ appears at most $g-1$ times.  Theorem \ref{thm1} recovers Corollary 2.1 in \cite{Cill-2010}.

If $A \subset [n]$ is a Sidon set, then for any $g \geq 2$, $A$ is a $C_2[g]$-set.  There are Sidon sets
$A \subset [n]$ with $|A| = (1  + o(1))n^{1/2}$ thus the exponent of (\ref{eq:BetterUpperBound-for-F_h(g,n)}) is correct when $h = 2$.  Motivated by constructions in extremal graph theory, we can show that (\ref{eq:BetterUpperBound-for-F_h(g,n)}) is correct for other values of $h$.

\begin{theorem}\label{sphere}
Let $p$ be an odd prime and $\alpha \in \mathbb{F}_p$ be chosen  to be a quadratic non-residue if $p \equiv 1 (\textup{mod}~4)$, and a nonzero quadratic residue otherwise.  The set
\[
A = \{ (x_1,x_2,x_3) \in \mathbb{F}_{p}^{3} : x_1^2 + x_2^2 + x_3^2 = \alpha \}
\]
is a $C_3[3]$-set in the group $\mathbb{F}_p^3$ with $|A| \geq p^2 - p$.
\end{theorem}

\begin{corollary}\label{c33 lb}
For any integer $n \geq 1$, there is a $C_3[3]$-set $A \subset [n]$ with
\[
|A| \geq \left( 4^{-2/3} + o(1) \right) n^{2/3}.
\]
\end{corollary}

By (\ref{group bound}), Theorem~\ref{sphere} is asymptotically best possible.  It is an open problem to determine the maximum size of a $C_3[3]$-set in $[n]$.

Proposition~\ref{prop1} suggests that the methods used to construct $K_{g,h}$-free graphs may be used to construct $C_h[g]$-sets.  Using the norm graphs of Koll\'{a}r, R\'{o}nyai, and Szab\'{o} \cite{krs}, we construct $C_h [ h! + 1]$-sets $A \subset [n]$ with $|A| \geq c_h n^{1 - 1/h}$ for each $h \geq 2$.

\begin{theorem}\label{nq lb}
Let $h \geq 2$ be an integer.  For any integer $n$, there is a $C_h[h!+1]$-set $A \subset [n]$ with
\[
|A| = (1 + o(1)) \left(  \frac{n}{ 2^{h-1} } \right)^{1 - 1/h}.
\]
\end{theorem}

Using the probabilistic method we can construct sets that are almost $C_h[g]$ for all $g \geq h \geq 2$.   A set $A \subset \Gamma$ is a \emph{weak $C_h[g]$-set} if given any set $X \subset \Gamma$ with $|X| = k$, and any set $\{ k_1 , \dots , k_g \} \subset \Gamma$ such that $X + k_1 , \dots , X + k_g$ are all pairwise disjoint, at least one of the translates $X+ k_i$ is not contained in $A$.
Erd\H os and Harzheim used the probabilistic method to construct such sets.  Here we do the same but obtain a better lower bound.

\begin{theorem} \label{thm2}
For any integers $g \geq h \geq 2$, there exists a weak-$C_h[g]$-set $A \subset [n]$ such that
$$|A| \geq \frac{1}{8} n^{\left( 1- \frac{1}{h}\right)\left( 1- \frac{1}{g} \right)\left(1+\frac{1}{hg-1} \right)}
.$$
\end{theorem}
It should be noted that for $h$ fixed, Theorem \ref{thm2} gives $|A|\geq  n^{1-\frac 1h-\epsilon}$ for $g$ sufficiently large, being a lower bound close to the exponent given in Theorem~\ref{thm1}.

Erd\H{o}s and Harzheim also proved that for any infinite $C_h[g]$-sequence $A \subset \mathbb{N}$,
\begin{equation*} \label{eq:UpperEstimate-C_h[g]-sequence}
\liminf_{n \to \infty} \frac{A(n)}{n^{1-1/h}} = 0.
\end{equation*}
Here $A(n) = | A \cap \{1,2, \dots , n\} |$.
We refine this result as follows.
\begin{theorem} \label{thm3}
If $A$ is an infinite $C_h[g]$-sequence with $g \geq h \geq 2$, then
\begin{equation*}
\liminf_{n \to \infty} \frac{A(n)(\log n)^{1/h}}{n^{1 - 1/h}} =O(1),
\end{equation*}
where the implicit constant depends only on $g$ and $h$.
\end{theorem}
Theorem \ref{thm3} was proved by Erd\H os \cite{Erd-1957} when $h=g=2$.

The rest of the paper is organized as follows.  In Section 2 we prove Theorem~\ref{thm1}.   We shall prove
Theorem \ref{sphere}, Corollary~\ref{c33 lb}, and Theorem \ref{nq lb} in Section 3.  Theorem~\ref{thm2} is proved in Section 4 and Theorem~\ref{thm3} is proved in Section 5. We conclude will some open problems.


\section{Proof of Theorem~\ref{thm1}}

We will use an inequality due to Cilleruelo and Tenenbaum~\cite{CillTen-2007}.

\begin{theorem}[Overlapping Theorem \cite{CillTen-2007}] \label{thm:overlapping-thm}
Let $(\Omega, \mathcal{A}, \mathbb{P})$ be a probability space and let $\{ E_j \}_{j=1}^k$ denote a family of events. For $m \geq 1$, let
\[
\sigma_m := \sum_{1 \le j_1 < \cdots < j_m \le k} \Pr( E_{j_1} \cap \cdots \cap E_{j_m}).
\]
Then for any $m \geq 1$,
\[
\sigma_m \ge \binom{\sigma_1}{m}=\frac{\sigma_1(\sigma_1-1)\cdots (\sigma_1-(m-1))}{m!}.
\]
\end{theorem}

\begin{proof}[Proof of Theorem~\ref{thm1}] Let $A \subset [n]$ be a $C_h[g]$-set  and let $B$ be any subset of $[n]$ with size at least $h$.
Let $\Y$ be a random variable with range the positive integers and law
\[
    \Pr(\Y=m) =
	\begin{cases}
 	\dfrac{1}{|A+B|} & \text{if }  m \in A+B, \\
 	0 & \text{otherwise. }
	\end{cases}
\]
For every $b \in B$ we define the event $E_b = \{ \omega \in \Omega \colon \Y(\omega) \in  A + b \},$ that has probability
$\Pr(E_b)=\sum_{a\in A}\Pr(\Y=a+b)= |A|/|A+B|$. We also write
\[
\sigma_m := \sum_{\{b_{1}, \cdots, b_{m}\} \in \binom{B}{m}} \Pr( E_{b_1} \cap \cdots \cap E_{b_m} ), \quad (m \ge 1).
\]
In particular
$$
 \sigma_1 = \frac{|A||B|}{|A+B|}.
$$
Let  $b_1 > \cdots > b_h$ be $h$ fixed elements of $B$. We can write
\begin{align*}
 \Pr(E_{b_1} \cap \cdots \cap E_{b_{h}})  & = \sum_{ \{a_1, \cdots, a_h  \} \in \binom{A}{h} } \Pr(\Y=a_1+b_1=a_2+b_2 \cdots =a_h+b_h) \\
 							& = \sum_{ a_1 + \{0,b_1-b_2, b_1-b_3, \cdots, b_1-b_h \}  \in \binom{A}{h}} \frac{1}{|A+B|},
\end{align*}
the sum extending to all $a_1\in A$ such that
$ a_1 + \{0,b_1-b_2, b_1-b_3, \cdots, b_1-b_h \}\subset A.$ These are congruent $h$-subsets of the  $C_{h}[g]$-set $A$, thus
\begin{align*} 							
							 \Pr(E_{b_1} \cap \cdots \cap E_{b_{h}})
							& \le \dfrac{g-1}{|A+B|}.
\end{align*}
Now we use Theorem \ref{thm:overlapping-thm}  to obtain
$$
	\binom{|B|}{h} \; \dfrac{g-1}{|A+B|} \ge \sigma_{h} \ge  \frac{\sigma_1 (\sigma_1-1) \cdots  (\sigma_1-h+1)}{h!}\ge \frac{\sigma_1}{h!}(\sigma_1-(h-1))^{h-1},
$$
and so
$$
 \frac{|B|^h}{h!} \; \dfrac{(g-1)}{|A+B|} \ge \frac{|A||B|}{h! \, |A+B|} \left(\frac{|A||B|}{|A+B|}-(h-1) \right)^{h-1},
$$
which implies
\[
    |A|^{h/(h-1)} \le | A + B | \left( (g-1)^{1/(h-1)} + \frac{(h-1)|A|^{1/(h-1)}}{|B|} \right).
\]
If we choose $B = [l]$, by the last inequality we have
\begin{equation} \label{eq:MejorCota-de-F_h(g,n)_1}
|A|^{h/(h-1)} \le (n+\ell) \left( (g-1)^{1/(h-1)} + \frac{(h-1)|A|^{1/(h-1)}}{\ell+1} \right).
\end{equation}
We first take  $\ell=n$ and use  $|A|\le n$ in the right side, getting
$$|A|^{h/(h-1)}=O(n) \implies |A|^{1/(h-1)}=O(n^{1/h}).$$
Inserting this in the second member of \eqref{eq:MejorCota-de-F_h(g,n)_1} we obtain
\[
|A|^{h/(h-1)} \le (g-1)^{1/(h-1)} n + O( \ell) \ + O\left (\frac{ n^{1+1/h}}{\ell+1}\right ) + O( n^{1/h}).
\]
To minimize this last upper bound we choose $\ell \asymp n^{1/2+1/2h}$. Then we can write
\[
|A|^{h/(h-1)} \le (g-1)^{1/(h-1)} n + O \left( n^{1/2+1/2h} \right) = (g-1)^{1/(h-1)} n \left( 1 + O\left( n^{1/2h-1/2} \right)\right),
\]
which yields
\[
|A| \le (g-1)^{1/h}  n^{1-1/h} \left( 1 + O\left( n^{1/2h-1/2} \right)\right)^{1 -1/h} = (g-1)^{1/h} n^{1-1/h} + O\left( n^{1/2 - 1/2h)} \right), \qedhere
\]
as we claimed.
\end{proof}


\section{Proof of Theorem \ref{sphere}, Corollary \ref{c33 lb}, and Theorem \ref{nq lb}}

\begin{proof} [Proof of Theorem \ref{sphere}]  Recall that we choose $\alpha \in \fp$  as  a quadratic non-residue when $p \equiv 1 (\textup{mod}~4)$ and a nonzero quadratic residue otherwise. Let $G=(V,E)$ be a graph
with vertex set $V = \fp^3$. For $x=(x_1,x_2,x_3)$ and $y=(y_1,y_2,y_3)$ we have $(x,y) \in E(G)$ if and only if
\[
\sum_{i=1}^3 (x_i-y_i)^2=\alpha.
\]
The graph $G$ is $K_{3,3}$-free as shown by Brown \cite{b}.  Define
\[
S(\alpha)=\{(x_1,x_2,x_3) \in \fp^3 \colon x_1^2+x_2^2+x_3^2=\alpha\}.
\]
 Let $X=\{x,y,z\} \subset \fp^3$. Suppose $X+a \subset S(\alpha)$ for some $a \in \fp^3$. We first show  ${-a} \not \in \{x,y,z\}$. If $x=-a$ then we get  $0=\alpha$ as $x+a \in S(\alpha)$.  This is a contradiction since we have chosen $\alpha$ so that
$\alpha \neq 0$.  Therefore $-a \neq x$ and similarly, $-a \neq y$ and $-a \neq z$.
By definition, $(x,-a), (y,-a),$ and $(z,-a)$ are three edges in $G$, which tell us that $a$ is a common neighbor of $x,y,$  and $z$.   Assume that there are three translates $X+a, X+b, X+c$  contained in $S(\alpha)$ for distinct $a,b,c \in \fp^3$. We have $\{x,y,z\} \cap \{a,b,c\}=\emptyset$, and so $L=\{x,y,z\}$ and $R=\{a,b,c\}$ form a $K_{3,3}$ in $G$.
However $G$ is $K_{3,3}$-free, a contradiction.  Thus there are at most two elements $a,b \in \fp^3$ such that the translates $X+a$ and  $X+b$ are contained in $S(\alpha)$. This holds  for every $X \subset \fp^3$ with $|X|=3$. We proved Theorem \ref{sphere}.
\end{proof}

Next we prove Corollary~\ref{c33 lb} and Theorem~\ref{nq lb}.  Both results rely on the following lemma.

\begin{lemma}\label{chg lemma 2}
Let $p$ be a prime and $d \geq 1$ be an integer.
Define $\phi : \fp^d \rightarrow \ZZ$ by
\[
\phi ( ( x_1 , \dots , x_d ) ) = x_1 + 2p x_2 + (2p)^2 x_3+ \dots + (2p)^{d-1}x_d.
\]
where $0 \leq x_i \leq p -1$.  The map $\phi$ is 1-to-1 and furthermore, for any $x,y,z,t \in \fp^d$, we have $x + y = z + t$ if and only if $\phi (x) + \phi (y) = \phi (z) + \phi (t)$.
\end{lemma}

The proof of Lemma~\ref{chg lemma 2} is not difficult.  In the language of additive combinatorics, the map $\phi$ is a Frieman isomorphism of order 2 (see \cite{tv}, Chapter 5, Section 3).

\begin{proof}[Proof of Corollary \ref{c33 lb}]
Let $n$ be a large integer.  Choose an odd prime $p$ with $4p^3 \leq n$ and $p$ as large as possible.
Let $S \subset \mathbb{F}_{p}^3$ be a $C_3[3]$-set in $\mathbb{F}_p^3$
with $|S| \geq p^2 - p$ guaranteed by Theorem \ref{sphere}.  Consider $A = \phi ( S)$ where
$\phi : \fp^3 \rightarrow \mathbb{Z}$ is the map
\[
\phi ( (x_1 , x_2 , x_3 ) ) = x_1 + 2p x_2 + 4p^2 x_3
\]
Here $x_i$ is chosen so that $0 \leq x_i \leq p -1$.  By Lemma~\ref{chg lemma 2}, $A$ is a $C_3[3]$-set.  If $a \in A$, then $a \leq (p-1)(1 + 2p + 4p^2) \leq 4p^3 \leq n$ so $A \subset [n]$.  Since $\phi$ is 1-to-1,
$|A| \geq p^2 - p$.  For large enough $n$, there is always a prime between $(n/4)^{1/3} - (n/4)^{ \theta / 3}$ and
$( n /4)^{ 1/3}$ for some $\theta < 1$.  The results of \cite{bhp} show that one can take $\theta = 0.525$.  Therefore,
$|A| \geq ( n /4)^{2/3} - O ( n^{ \frac{\theta  + 1}{3}} ) = (1 + o(1))(n/4)^{2/3}$.
\end{proof}

\begin{proof}[Proof of Theorem~\ref{nq lb}]

Let $q$ be a prime power and $h \geq 2$ be an integer.  Let $N: \mathbb{F}_{q^h} \rightarrow \mathbb{F}$ be the norm map defined by
\[
N(x) = x^{ 1 + q + q^2 + \dots + q^{h- 1} }.
\]
Let $A = \{ x \in \mathbb{F}_{q^h} : N(x)  = 1 \}$.  The norm map $N$ is a group homomorphism that maps
$\mathbb{F}_{q^{h}}^*$ onto $\mathbb{F}_{q}^{*}$.  This implies $\frac{q^h - 1}{ |A| } = q  -1$
so $|A| = \frac{q^h - 1}{q-1}$.  We now show that $A$ is a $C_h [h! + 1]$-set in the group $\mathbb{F}_{q^h}$.

Suppose $X = \{ x_1 , \dots , x_h \} \subset \mathbb{F}_{q^h}$.  It follows from Theorem 3.3 of \cite{krs} that there are at most $h!$ elements $k \in \mathbb{F}_{q^h}$ such that
\[
N( k + x_i) =1
\]
for all $1 \leq i \leq h$.  Therefore, given any set $\{k_1 , \dots , k_{h! +1} \} \subset \mathbb{F}_{q^h}$,
at least one of the translates $X + k_i$ is not contained in $A$.

Let $\psi : \mathbb{F}_{q^h} \rightarrow \mathbb{Z}_q^h$ be a group isomorphism mapping the additive group
$\mathbb{F}_{q^h}$ onto the direct product $\mathbb{Z}_q^h$.
Let $\phi : \mathbb{Z}_{q}^{h} \rightarrow \mathbb{Z}$ be the map
\[
\phi (x_1 , \dots , x_h ) = x_1 + (2q)x_2 + \dots + (2q)^{h-1} x_h
\]
where $0 \leq x_i \leq q - 1$.  By Lemma~\ref{chg lemma 2}, $A' := \phi ( \psi (A) )$ is a
$C_h [h! + 1]$-set.  The set $A'$ has $\frac{q^h - 1}{q-1}$ elements and is contained in
the set $[ 2^{h-1} q^h ]$.  By the same argument used to prove Corollary \ref{c33 lb}, we can choose a prime power $q$ given a large enough integer $n$ to obtain a $C_h[h!+1]$-set in $[n]$ with size $(1+o(1))\left( \tfrac{n}{2^{h-1}}\right)^{1-1/h}$.
\end{proof}


\section{Proof of Theorem~\ref{thm2}}
The proof in this section uses the probabilistic method combined with the deletion technique.
These ideas have appeared before in the literature, see for example \cite[\S 3]{AloSpe-2000}, \cite{SpeTet-1995},  and \cite{Cill-2010a}.

We say that
$m \in S$ is $(h,g)$-bad (for $S$) if there exist $m_1 < \cdots < m_{g-1}$, with $m_i < m$, and there exist $\ell_1<\ell_2 < \cdots < \ell_{h-1}$ such that the sums $\{m_1,\cdots, m_{g-1},m \} + \{0,\ell_1, \cdots, \ell_{h-1}\} $ are $gh$ distinct elements of $S$.

We define $S_{bad}$ the set of $(h,g)$-bad elements for $S$.
It is clear that for any set $S$, the set $$S_{C_h[g]}=S\setminus{S_{bad}},$$ is weak-$C_h[g]$-set with cardinality $|S_{C_h[g]}|=|S|-|S_{bad}|.$

Define $p$ as the number such that $2pn=n^{g+h-1}(2p)^{hg}$. It is straightforward to check that \begin{equation}\label{np}np=\frac 12n^{\left( 1- \frac{1}{h}\right)\left( 1- \frac{1}{g} \right)\left(1+\frac{1}{hg-1} \right)}.\end{equation}
We will prove that except for finitely many $n$ there exist a set $S \subset [n]$ such that
\begin{equation}\label{ine}
|S|\ge\frac{np}2\quad \text{ and }\quad |S_{bad}|\le\frac{np}4.
\end{equation}
Note that for such a set we have
$$|S_{C_h[g]}|=|S|-|S_{bad}|>\frac{np}4=\frac{1 }{8} n^{\left( 1- \frac{1}{h}\right)\left( 1- \frac{1}{g} \right)\left(1+\frac{1}{hg-1} \right)},$$ for all sufficiently large $n$ and $A=S_{C_h[g]}$ satisfies the conditions of Theorem \ref{thm2}.

Indeed we will prove that  with probability at least $1/4$, a random set $S$  in $[n]$ satisfies \eqref{ine} if each element in $[n]$ is independently chosen to be in $S$ with probability $p$.

Next we obtain estimates for the random variables $|S|$ and $|S_{bad}|$.

If $m$ is $(h,g)$-bad then the $gh$ sums $\{m_1,\cdots, m_{g-1}, m \} + \{0,\ell_1, \cdots, \ell_{h-1}\}$ are all distinct elements of $S$ and so
$$\Pr(\{m_1,\dots,m_{g-1},m\}+\{0,\ell_1,\dots,\ell_{h-1}\}\subset S)= p^{gh}.$$
Hence
\begin{align*}
&\Pr(m \text{ is $(h,g)$-bad} ) \le
\sum_{ \substack{1 \le m_1 < \cdots < m_{g-1} < m \\1 \le \ell_1 < \cdots < \ell_{h-1} \le n } }  p^{gh}\le \binom m{g-1}\binom n{h-1}p^{gh}<n^{g+h-2}p^{gh}
\end{align*}
which implies
\[
 \E(|S_{bad}|) \le \sum_{1 \le m \le n} \Pr(m \text{ is $g$-bad} ) \le n^{g+h-1}p^{gh} .
\]
On the one hand by Markov's inequality we have
\begin{align}\label{Markov}
\Pr\left(|S_{bad}|>\frac{np}4 \right)	&=   \Pr\left(|S_{bad}|>\frac{n^{g+h-1}(2p)^{gh}}{8} \right)\\
										&=   \Pr(|S_{bad}|>2^{gh-3}n^{g+h-1}p^{gh})\nonumber\\
										&\le \Pr\left(|S_{bad}|>2\E(|S_{bad}|)\right)<1/2.\nonumber
\end{align}
On the other hand, using that $\E(|S|)=np$ and $\text{Var}(|S|)=np(1-p)$ and applying Chebychev's inequality we have
\begin{eqnarray}\label{Che}
\Pr\left(|S|<\frac{np}2 \right)&=&\Pr\left(|S|<\frac{\E(|S|)}2\right)<\Pr\left(|S-\E(|S|)|>\frac{\E(|S|)}2 \right)\\ &<&\frac{4\text{Var}(|S|)}{(\E(|S|))^2}=\frac{4np(1-p)}{(pn)^2}<\frac{4}{pn}<\frac 14,\nonumber
\end{eqnarray} except for finitely many $n$.
By \eqref{Markov} and \eqref{Che} we have
$$ \Pr(|S|\ge np/2 \text{ and } |S_{bad}|\le np/4) \ge 1-\left (1/2+1/4\right )\ge 1/4, $$
as we wanted.

\section{Proof of Theorem \ref{thm3}}
In order to simplify notation, when $f(n)=O(g(n))$ we write $f(n) \ll g(n)$ or $g(n) \gg f(n)$ through this section.

\begin{proof}[Proof of Theorem \ref{thm3}] Let $A$ be an infinite $C_h[g]$-sequence.
For a positive integer $N$, let $[N^2]$ denote all the positive integers less or equal to $N^2.$
We divide $[N^2]$ into equally sized intervals $$I_\nu := [(\nu-1)N, \nu N], \; \nu=1, \cdots, N.$$
Let $\mathcal{C}$ denote the collection of all $h$-subsets of $[N^2]$ that are included in one of the intervals $I_\nu$:
\[
 \mathcal{C} := \left \lbrace C \in \binom{[N^2]}{h} \colon C \subset I_\nu \text{ for some } \nu \right \rbrace.
\]
We say that the sets in the collection $\mathcal{C}$ are \lq\lq small\rq\rq \ as their diameter is at most $N$.
We classify the elements of $\mathcal{C}$ so that each class groups all the sets that are pairwise congruent.
Each class $\alpha $ contains a set $C_\alpha$ that contains $1$,
and the remaining $h-1$ elements of $C_\alpha$ can be chosen in $\binom{N-1}{h-1}$ different ways;
each of the choices determines a class different from the others. Then the number of classes is
\[
  \binom{N-1}{h-1}.
\]
Let $A_\nu$ denote the size of $A \, \cap \, I_\nu$, we have $A_\nu = A(\nu N) - A((\nu-1)N)$, where $A(x) := |\{ a \in A \colon a \le x \}|$ is the counting function of the sequence.

One the one hand as $A$ is a $C_h[g]$-sequence then in every class of $\mathcal{C}$ there are at most $g-1$ subsets of $A$. 
Hence we have the following upper bound for the total number of \lq\lq small\rq\rq \ subsets of $A$ that belong to $\mathcal{C}$
\begin{equation*} \label{eq:liminf-ll-1-1}
\sum_{\nu=1}^{N} \binom{A_\nu}{h} \le \binom{N-1}{h-1} (g-1) \ll N^{h-1} \quad (N \to \infty),
\end{equation*}
Now we prove by induction in $h$ that
\begin{equation} \label{eq:liminf-ll-1-3}
\sum_{\nu=1}^{N} A_\nu^h \ll N^{h-1} \quad (N \to \infty).
\end{equation}
For $h=2$ we know by Theorem \ref{eq:BetterUpperBound-for-F_h(g,n)} that $A(N^2)\ll N$, so
$$\sum_{\nu=1}^{N} A_\nu^2 = 2 \sum_{\nu=1}^{N} \binom{A_\nu}{2} + \sum_{\nu=1}^{N} A_\nu \ll N + A(N^2) \ll N.$$
If \eqref{eq:liminf-ll-1-3} holds for all exponents up to $h-1$, then
\begin{equation*}
\sum_{\nu=1}^{N} A_\nu^h = h! \sum_{\nu=1}^{N} \binom{A_\nu}{h}+ O\left( \sum_{\nu=1}^{N} A_\nu^{h-1} \right) \ll N^{h-1} +N^{h-2}, \qquad (N \to \infty),
\end{equation*}
thus it also holds for $h$. Using \eqref{eq:liminf-ll-1-3} and  H\"{o}lder inequality we can write
\begin{align}
\sum_{\nu=1}^N A_\nu \left(\frac{1}{\nu}\right)^{1-1/h} & \le \left( \sum_{\nu=1}^N A_\nu^h \right)^{1/h} \, \left( \sum_{\nu=1}^N \frac{1}{\nu} \right)^{1-1/h} \notag \\
		&\ll \left( N \log N \right)^{1-1/h}, \qquad (N \to \infty) \label{eq:liminf-ll-1-4}.
\end{align}

\medskip
On the other hand as $\sum_{\nu \le t} A_\nu=A(tN)$ and summing by parts
\[
\sum_{\nu=1}^N A_\nu \left(\frac{1}{\nu}\right)^{1-1/h} = \frac{A(N^2)}{N^{1-1/h}}+\int_1^N \frac{A(tN)}{t^{2-1/h}} \d t.
\]
In this sum the first summand is bounded by Theorem \ref{thm1} as follows
\begin{equation}
\frac{A(N^2)}{N^{1-1/h}} \ll N^{2(1-1/h)}/N^{(1-1/h)}=N^{(1-1/h)},\label{eq:last}
\end{equation}
and as consequence we shall prove next that the second summand is the main term in the sum.
Let us write $$\tau(m):=\inf_{n \ge m} \frac{A(n)(\log n)^{1/h}}{n^{1-1/h}}.$$
For $N \geq m$ and $t \geq 1$ we have
\[
A(tN)=\frac{A(tN)(\log(tN))^{1/h}(tN)^{1-1/h}}{(tN)^{1-1/h} (\log(tN))^{1/h}} \geq \tau(m) \frac{t^{1-1/h}N^{1-1/h}}{(\log N)^{1/h}}.
\]
Thus for $N \ge m$ we have
\[
\int_1^N \frac{A(tN)}{t^{2-1/h}} \d t \gg \frac{\tau(m)N^{1-1/h}}{(\log N)^{1/h}} \int_1^N \frac{1}{t} \d t \gg \tau(m) (N \log N)^{1-1/h},
\]
and so by \eqref{eq:last}
\begin{equation*}
\sum_{\nu=1}^N A_\nu \left(\frac{1}{\nu}\right)^{1-1/h} \gg \int_1^N \frac{A(tN)}{t^{2-1/h}} \d t \gg \tau(m) (N \log N)^{1-1/h}.
\end{equation*}
Inserting $\eqref{eq:liminf-ll-1-4}$ we have $\lim_{m \to \infty} \tau(m) \ll 1,$
that is what we wanted to prove.
\end{proof}

\section{Open problems}

In this final section we mention several open problems.

\indent

\textbf{Problem 1:}  Determine the maximum size of a $C_3[3]$-set contained in $[n]$.

Our results show that if $A \subset [n]$ is a $C_3[3]$-set of maximum size then
\[
(4^{-2/3} + o(1)) n^{2/3} \leq |A| \leq (2^{1/3} + o(1)) n^{2/3}.
\]
It seems likely that both of these bounds can be improved.  Perhaps the correct answer is $(1 + o(1))n^{2/3}$.

\textbf{Problem 2:} Remove the condition \emph{weak} in Theorem \ref{thm2}.

A much harder problem is the following.

\textbf{Problem 3:}  Construct $C_h[g]$-sets in $[n]$ with the order $n^{1-\tfrac{1}{h}}$ for each $g \geq h \geq 3$.

For $g \geq h=3$ and $g \geq h!+1$, we constructed $C_h[g]$-sets in $[n]$ whose sizes matches the order given by Theorem \ref{thm1}. We believe for any other $g$ and $h$ the upper bound by Theorem \ref{thm1} gives the correct exponent.
We note that solving Problem 3 would imply $z(n,n,g,h) \geq C(g,h) n^{1-1/h}$ for some constant $C$  (see Proposition \ref{prop1}).

\textbf{Problem 4:} Construct an infinite $C_h[g]$-sequence  $A \in \mathbb N$ which has  counting function $A(n)\gg n^{\left( 1- \frac{1}{h}\right)\left( 1- \frac{1}{g} \right)\left(1+\frac{1}{hg-1} \right)+o(1)}$ for all $n$.

We have found technical difficulties to deal with Problems 2 and 4, which were suggested to the second author by Javier Cilleruelo.



%


%
%

\end{document}